\documentclass[10pt,francais]{smfart}

\usepackage[francais]{babel}

\usepackage{smfthm}

\newcommand{\C}{\mathbb C}
\renewcommand{\P}{\mathbb P}
\newcommand{\N}{\mathbb N}

\newtheorem{theorem}{Th\'eor\`eme}[section]
\newtheorem{lemma}[theorem]{Lemme}
\newtheorem{proposition}[theorem]{Proposition}
\newtheorem{corollary}[theorem]{Corollaire}
\newtheorem{rappel}[theorem]{}
\newtheorem{problem}[theorem]{Probl\`eme}

\theoremstyle{definition}
\newtheorem{definition}[theorem]{D\'efinition}
\newtheorem{remark}[theorem]{Remarque}

\newtheorem{notation}[theorem]{Notation}
\newtheorem{question}[theorem]{Question}
\newtheorem{example}[theorem]{Exemple}

\newcommand{\cal}{\mathcal}
\newcommand{\wh}{\widehat}
\newcommand{\lan}{\langle}
\newcommand{\ran}{\rangle}

\newcommand{\bi}{\begin{itemize}}
\newcommand{\ei}{\end{itemize}}
\newcommand{\be}{\begin{enumerate}}
\newcommand{\ee}{\end{enumerate}}

\newcommand{\bpf}{\begin{proof}}
\newcommand{\epf}{\end{proof}}
\newcommand{\bt}{\begin{theorem}}
\newcommand{\et}{\end{theorem}}
\newcommand{\brap}{\begin{rappel}}
\newcommand{\erap}{\end{rappel}}
\newcommand{\bnt}{\begin{notation}}
\newcommand{\ent}{\end{notation}}
\newcommand{\bd}{\begin{definition}}
\newcommand{\ed}{\end{definition}}
\newcommand{\ble}{\begin{lemma}}
\newcommand{\ele}{\end{lemma}}
\newcommand{\bpr}{\begin{proposition}}
\newcommand{\epr}{\end{proposition}}
\newcommand{\bre}{\begin{remark}}
\newcommand{\ere}{\end{remark}}
\newcommand{\bco}{\begin{corollary}}
\newcommand{\eco}{\end{corollary}}
\newcommand{\beq}{\begin{equation}}
\newcommand{\eeq}{\end{equation}}
\newcommand{\bq}{\begin{question}}
\newcommand{\eq}{\end{question}}
\newcommand{\bp}{\begin{problem}}
\newcommand{\ep}{\end{problem}}
\newcommand{\beqn}{\begin{eqnarray*}}
\newcommand{\eeqn}{\end{eqnarray*}}
\newcommand{\bex}{\begin{example}}
\newcommand{\eex}{\end{example}}
\newcommand{\sct}{\section}
\newcommand{\ssct}{\subsection}
\newcommand{\sk}{\smallskip}
\newcommand{\bk}{\bigskip}
\newcommand{\nk}{\noindent}

\newcommand{\fr}{\frac}

\newcommand{\bck}{\backslash}

\renewcommand{\a}{\alpha}
\renewcommand{\b}{\beta}

\renewcommand{\O}{\Omega}
\renewcommand{\k}{\kappa}

\newcommand{\s}{{\sigma}}
\newcommand{\Si}{{\Sigma}}

\begin{document}

\begin{abstract}
Soit $X$ un germe de vari\'et\'e complexe lisse en $x\in \P^N$,
\`a osculation $q$-r\'eguli\`ere et tel que, pour  toute direction 
$v\in \P T_xX$, il existe une courbe rationnelle normale 
de degr\'e $q$ tangente en $x$ \`a $v$ et localement contenue dans $X$. 
On montre que $X$ est une vari\'et\'e de Veronese 
d'ordre $q$. On retrouve en particulier une caract\'erisation des vari\'et\'es 
de Veronese, due \`a Bompiani \cite{Bo}, dont on commente les r\'esultats.
\end{abstract}

\title[Une caract\'erisation des vari\'etes de Veronese]{une nouvelle caract\'erisation\\
des vari\'et\'es de Veronese}
\author{Jean-Marie Tr\'epreau}

\address{U.P.M.C., 
UMR 7586,  bureau 26-16-5-21 ; 
4, place Jussieu, 75005 Paris.
Adresse \'electronique : trepreau@math.jussieu.fr}

\maketitle

\sct{Introduction}

\ssct{\'Enonc\'es}

Soit $n\in \N^*$ et $q\in \N^*$. Une sous-vari\'et\'e $X$ de dimension $n$ d'un espace projectif 
est une {\em vari\'et\'e de Veronese} d'ordre $q$ si, dans l'espace projectif qu'elle engendre,
$X$ est l'image de $\P^n$ par un plongement associ\'e au syst\`eme lin\'eaire $\left| \cal{O}_{\P^n}(q) \right|$.

\sk
L'image d'une droite de  $\P^n$ par ce plongement est une vari\'et\'e 
de Veronese de dimension 1 et d'ordre $q$, autrement dit une 
{\em courbe rationnelle normale} de degr\'e $q$. Pour toute paire de points distincts 
d'une vari\'et\'e de Veronese $X$
d'ordre $q$, il existe donc une courbe rationnelle normale de degr\'e $q$,
contenue dans $X$ et qui passe par ces points. D'autre part, la vari\'et\'e $X$
engendre un espace de dimension $N = { n+q \choose n} - 1$, la dimension de
l'espace des polyn\^omes homog\`enes de degr\'e $q$ en $n+1$ variables.

\sk
L'espace osculateur $X_x(q)$ \`a l'ordre $q$ \`a un germe de vari\'et\'e lisse 
$X\subset \P^N$ en $x\in X$ est le sous-espace projectif de $\P^N$, passant par $x$
et de direction l'espace engendr\'e au point $x$ par les d\'eriv\'ees d'ordre $\leq q$
des germes de courbes param\'etr\'ees trac\'ees sur $X$.
Nous disons que {\em la vari\'et\'e $X$ est $q$-r\'eguli\`ere}  en $x$
si cet espace est de dimension maximale, c'et-\`a-dire de dimension ${n+q \choose q}-1$,
o\`u $n$ est la dimension de $X$.

\sk
Dans \cite{Bo}, E. Bompiani obtenait une caract\'erisation g\'eom\'etrico-diff\'erentielle 
des vari\'et\'es de Veronese qui contient, comme cas 
particulier, l'\'enonc\'e suivant. Nous reviendrons dans la Section 4 sur l'article de Bompiani.
\bt
\label{Th0}
Soit $X$ un germe de vari\'et\'e lisse et $q$-r\'eguli\`ere. 
Si pour toute paire  $(x,y)\in X^2$
g\'en\'erique, il existe une courbe rationnelle normale de degr\'e $q$
localement contenue dans $X$ et passant par $x$ et $y$, alors $X$
est le germe d'une vari\'et\'e de Veronese d'ordre $q$.
\et
Nous obtiendrons, par un calcul formel, le r\'esultat suivant :
\bt
\label{Th1}
Soit $X$ un germe de vari\'et\'e lisse et $q$-r\'eguli\`ere 
en $x\in \P^N$. Si pour tout $v\in \P T_xX$,
il existe un germe de courbe rationnelle
normale de degr\'e $q$ contenu dans $X$ et tangent \`a $v$ en $x$,
$\,X$ est le germe d'une vari\'et\'e de Veronese d'ordre $q$. 
\et
\bre
La d\'emonstration montre qu'il suffit que l'ensemble des
$v$ qui ont la propri\'et\'e de l'\'enonc\'e 
ne soit pas contenu dans une hypersurface alg\'ebrique de $\P T_xX$.
\ere
La condition de $q$-r\'egularit\'e de l'\'enonc\'e 
est n\'ecessaire. On v\'erifie en effet que 
la projection d'une vari\'et\'e de Veronese de dimension $\geq 2$ 
depuis un point g\'en\'erique de
l'espace qu'elle engendre v\'erifie l'hypoth\`ese principale de l'\'enonc\'e
et n'est pas $q$-r\'eguli\`ere.

\sk 
On peut, dans la caract\'erisation de Bompiani, se donner le point $x$ :
\bt
\label{Th1bis}
Soit $X$ un germe de vari\'et\'e lisse et $q$-r\'egul\`ere en $x\in \P^N$. Si 
pour tout  $y\in X$ g\'en\'erique, il existe une courbe rationnelle
normale de degr\'e $q$ localement con\-tenue dans $X$ et passant
par $x$ et $y$, $\,X$ est le germe d'une vari\'et\'e de Veronese d'ordre~$q$.
\et

\ssct{Motivations et un probl\`eme ouvert}

Le th\'eor\`eme de Bompiani joue un r\^ole important
dans le probl\`eme de l'alg\'ebrisation des tissus de rang maximal,
sur lequel je travaille depuis plusieurs ann\'ees  avec Luc Pirio, voir Pirio-Tr\'epreau \cite{PT2}.

\sk
Nous avions en effet besoin, pour r\'esoudre le probl\`eme
principal, de r\'esoudre le probl\`eme auxiliaire de g\'eom\'etrie projective suivant,
o\`u $r\geq 1$, $n\geq 2$ et $q\geq n-1$  sont des entiers donn\'es.

\sk
{\em 
D\'eterminer toutes les vari\'et\'es projectives $X$ de dimension $r+1$ telles que :
\be

\item[1)] pour tout $n$-uplet $(x_1,\ldots,x_n)\in X^n$ g\'en\'erique, il existe 
une courbe rationnelle normale de degr\'e $q$ contenue dans $X$ et passant par
$x_1,\ldots,x_n$ ;

\item[2)] $X$ engendre un espace projectif dont la dimension est la plus grande 
possible, pour $r$, $n$ et $q$ donn\'es, compte tenu de la propri\'et\'e 
pr\'ec\'edente.

\ee}

\sk 
Le probl\`eme est r\'esolu si $q\neq 2n-3$, dans Pirio-Tr\'epreau \cite{PT1}. 
Si $n=2$ le r\'esultat  n'est autre que la caract\'erisation par Bompiani 
des vari\'et\'es de Veronese. 

\sk
Parce que, dans  l'article de Bompiani, je n'ai pas trouv\'e convaincantes 
les quelques lignes consacr\'ees \`a la d\'emonstration du Th\'eor\`eme \ref{Th0},
dans le doute, j'en ai cherch\'e une qui me satisfasse. J'ai tent\'e d'en trouver 
une d\'emonstration formelle et j'ai obtenu
le Th\'eor\`eme \ref{Th1}. C'est cette d\'emonstration,
\'el\'ementaire et un peu laborieuse, qui est pr\'esent\'ee 
ici. 

\sk
Le Th\'eor\`eme \ref{Th1} est bien s\^ur plus fort 
que l'\'enonc\'e de Bompiani mais celui-ci est suffisant
dans \cite{PT1},\cite{PT2}. 
Des g\'en\'eralisations d'une autre nature auraient des cons\'equences int\'eres\-santes
pour \'etendre certains des r\'esultats de \cite{PT2} \`a des tissus qui ne sont pas de rang
maximal. Consid\'erons \`a nouveau une vari\'et\'e projective $X$ de dimension $n$, telle 
que pour toute paire g\'en\'erique $(x,y)$ de points de $X$, il existe 
une courbe rationnelle normale de degr\'e $q$ contenue dans $X$ et qui 
passe par les points $x,y$. Elle engendre un espace de dimension 
au plus $N= {n+q \choose q}-1$ et,
si elle engendre un espace de dimension $N$, c'est une vari\'et\'e de Veronese 
d'ordre $q$. La question suivante semble ouverte et n'est pas gratuite.

\sk
{\em Quelle est la valeur minimale de l'entier $M_{n,q}$ qui assure qu'on ait :
si $X$ engendre un espace de dimension $\geq M_{n,q}$, alors $X$ est une  projection 
d'une vari\'et\'e de Veronese d'ordre $q$ ?}

\ssct{Plan de l'article}
L'essentiel consiste en un calcul formel.
Dans la Section 2, nous \'enon\c{c}ons le r\'esultat de ce calcul, 
le Lemme \ref{conclusion}, dont nous d\'eduisons le Th\'eor\`eme~\ref{Th1}
et une g\'en\'eralisation, le Th\'eor\`eme \ref{Th2}, voir aussi la Remarque \ref{comment}.
Le Th\'eor\`eme~\ref{Th1bis} est un corollaire g\'eom\'etrique du Th\'eor\`eme \ref{Th1}.
Le calcul formel est une r\'ecurrence qui fait l'objet de la Section 3.
Enfin, dans la Section 4, nous faisons quelques commentaires sur 
l'article de Bompiani.

\bk
\nk
{\bf \em Remerciements.--- }
L'auteur remercie Luc Pirio, qui a d\'ecouvert l'existence de l'article de
Bompiani et  compris l'importance du Th\'eor\`eme \ref{Th0}
pour l'\'etude du probl\`eme qui nous occupait.

\sct{Un lemme et ses cons\'equences}

\ssct{Notations}
\label{notations}

On travaillera pr\`es de $0$ dans $\C^N$, o\`u $N= {n+q\choose q}-1$. 
On note $t$ le point g\'en\'eral de $\C$ et $s=(s_1,\ldots,s_n)$ le point g\'en\'eral de $\C^n$.
Si $\a=(\a_1,\ldots,\a_n)\in \N^n$,
on note $|\a|=\a_1+\cdots+\a_n$ la longueur de $\a$ et $s^\a=s_1^{\a_1}\ldots s_n^{\a_n}$.

\sk
On indexe les composantes de $x\in \C^N$ par les multiindices 
$\a\in \N^n$ dont la longueur est comprise entre $1$ et $q$, soit $x=(x_\a)_{1\leq |\a|\leq q}$.
On dit que $x_\a$ est une composante de poids $|\a|$ de $x$.

\ssct{La propri\'et\'e (P)}
Les entiers $n\geq 2$ et $q\geq 2$ sont fix\'es
et $X$ est un germe de vari\'et\'e {\em lisse} de dimension $n$ en $0\in \C^N$.
\bd
\label{P}
On appelle {\em courbe distingu\'ee de $X$}
tout germe de courbe lisse en $0\in \C^N$, contenu dans l'intersection 
de $X$ et d'un espace de dimension $q$. 
On dit que $X$ a la {\em propri\'et\'e (P)} 
si le germe $X$ est $q$-r\'egulier en $0$  et si, 
pour toute direction  $v\in \P T_0X$, 
il existe une courbe distingu\'ee de $X$ tangente \`a $v$.
\ed
Dans la Section 3, on d\'emontrera le r\'esultat suivant :
\ble
\label{conclusion}
Si $X$  a la propri\'et\'e $(P)$, il existe une homographie $\phi\in \text{\rm Aut}\, \P^N$
avec $\phi(0)=0$, telle que $\phi(X)$ admette une param\'etrisation de la forme suivante :
\beq
\label{normal}
x_\a = s^\a,  \;\; |\a|=1 \; ; \; \;\; x_\a = s^\a + g_a(s) \;\; \text{avec} \;\; g_\a(s)= O(|s|^{q+3}),
\;\; 2\leq |\a|\leq q.
\eeq
Les courbes distingu\'ees de $X$ sont alors les images des germes de droite
en $0\in \C^n$ par la param\'etrisation.
\ele

\ssct{D\'emonstration du Th\'eor\`eme \ref{Th1}}

Montrons d'abord le lemme suivant, dans l'esprit, formel,  de cet article :
\ble
Si le germe en $0\in \C^q$ d'une courbe rationnelle normale de degr\'e $q$
est donn\'e par 
$$
x_j = x_1^j + g_j(x_1) \;\; \text{avec} \;\; g_j(x_1)  = O(x_1^{q+3}), \;\;\; j=2,\ldots,q,
$$
alors $g_j(x_1)\equiv 0$ pour $j=2,\ldots,q$.
\ele
\bpf
Le cas $q=1$ est trivial. On suppose $q\geq 2$. Par d\'efinition d'une courbe rationnelle
normale, le germe admet une param\'etrisation
$$
x_j(t) = X_j(t)/X_0(t), \qquad j=1,\ldots,q,
$$
o\`u $X_0(t),\ldots,X_q(t)$ sont des polyn\^omes de degr\'e $\leq q$.
Un changement affine de param\`etre $t$ permet de supposer $X_0(t)=1+O(t)$
et $X_1(t)=t + O(t^2)$.

\sk
Montrons par r\'ecurrence sur $r\geq 1$ qu'on a 
\beq
\label{hyp-r}
X_j(t) =   t^j +  O(t^{j+r}), \qquad j=0,\ldots,q.
\eeq
C'est vrai si $r=1$ compte tenu de la normalisation pr\'ec\'edente.
Puisque les $X_j(t)$ sont de degr\'es $\leq q$,
il s'agit de montrer que c'est vrai pour $r=q+1$.
Par hypoth\`ese :
$$
X_j(t)X_0(t)^{j-1} = X_1(t)^j + 0(t^{q+3})
$$
pour $2\leq j\leq q$. On se donne $r\in \{1,\ldots,q\}$, on suppose 
(\ref{hyp-r}) v\'erifi\'e et on \'ecrit
$$
X_j(t)=t^j + a_jt^{j+r} + O(t^{j+r+1}), \qquad j=0,\ldots,q.
$$
Pour $j+r\leq q+2$, on d\'eveloppe la relation pr\'ec\'edente \`a l'ordre $j+r$ :
$$
(t^j + a_jt^{j+r})(1+a_0t^r)^{j-1} = (t+a_1t^{r+1})^j + 0(t^{j+r+1}), \qquad j\geq 2, \;\; j+r\leq q+2.
$$
On obtient (le premier groupe d'\'equations parce que $X_j(t)$ est de degr\'e $\leq q$) :
\beqn
a_j          & = & 0,     \qquad \;\;\, j\geq 1,   \;\;\;  j+r\geq q+1    \\
a_j+(j-1)a_0 & = & ja_1,  \qquad        j\geq 2,   \;\;\;  j+r\leq q+2.
\eeqn
Si $r\leq q-1$, on peut choisir $j$ tel que $j+r=q+1$ puis tel que $j+r=q+2$ dans 
le deuxi\`eme groupe d'\'equations, ce qui donne $a_0=a_1=0$.
Si $r = q$, en choisissant $j=1$ dans le premier groupe puis $j=2$ dans le second,
on obtient encore $a_0=a_1=0$.

\sk
Dans tous les cas, on a $a_0=a_1=0$ et les \'equations donnent $a_j=0$
pour tout $j$. On obtient que (\ref{hyp-r}) est vrai \`a l'ordre $r+1$.
Par r\'ecurrence, le lemme est d\'emontr\'e.
\epf
Soit maintenant $X$ un germe de vari\'et\'e lisse et $q$-r\'eguli\`ere 
en $0\in \C^N$ tel que pour tout $v\in \P T_xX$,
il existe un germe de courbe rationnelle
normale de degr\'e $q$ contenu dans $X$ et tangent \`a $v$ en $x$.
On peut supposer $X$ donn\'e sous la forme r\'eduite du Lemme \ref{conclusion}.
Pour $\s=(\s_1,\ldots,\s_n)\in \C^n\bck\{0\}$, la courbe distingu\'ee 
$$
\xi_\a(t)  = \s^\a t \;\; \text{si}\;\; |\a|=1,
\qquad 
\xi_\a(t)  := \s^\a t^{|\a|} + g_\a(\s t) = \s^\a t^{|\a|} + O(t^{q+3}), \;\; 2\leq |\a| \leq q,
$$
est rationnelle normale  de degr\'e $q$. Pour $\a$ donn\'e, en appliquant le 
lemme pr\'ec\'edent aux $q$ composantes $\xi_\a(t)$ et $\xi_{(k,0,\ldots,0)}(t)$
avec $k\in \{1,\ldots,q\}\bck\{|\a|\}$, on obtient
$g_\a (\s t)\equiv 0$. Ainsi $X$ est donn\'e par 
$x_\a = s^\a$ pour $|\a|=1,\ldots,q$. \;C'est la vari\'et\'e de Veronese de dimension $n$
et d'ordre $q$ standard. 
Le Th\'eor\`eme \ref{Th1} est d\'emontr\'e.

\ssct{Germes qui ont la propri\'et\'e (P)}

On a le r\'esultat suivant :
\bt
\label{Th2}
Un germe de vari\'et\'e lisse en $0\in \C^N$ 
a la propri\'et\'e (P), voir la D\'efinition \ref{P},
si et seulement s'il est \'equivalent modulo $\text{\rm Aut}\, \P^N$
\`a un germe en $0$ d\'efini par une param\'etrisation de la forme,
avec $R_k(0)=0$ pour $k=1,\ldots,q$ : 
$$
x_\a = s^\a + s^\a R_{|\a|}(s), \qquad 1\leq |\a|\leq q.
$$
On peut imposer qu'on ait $R_1(s)\equiv 0$ et $s^\a R_{|\a|}(s)=O(|s|^{q+3})$ pour $|\a|\geq 2$.
\et
\bpf
Supposons d'abord le germe $X$ param\'etr\'e par
$$
x_\a = s^\a + s^\a R_{|\a|}(s), \qquad 1\leq |\a|\leq q,
$$
avec $R_k(0)=0$ pour $k=1,\ldots,q$. Il est $q$-r\'egulier en $0\in \C^N$. 
Pour $\s\in \C^n\bck\{0 \}$ fix\'e, l'image $\xi(t)$ du germe de droite 
d'\'equation $s=\s t$ est donn\'ee par 
$$
\xi_\a(t) = \s^\a t^{|\a|}  + \s^\a t^{|\a|}R_{|\a|}(\s t) = \s^\a S_{|\a|}(t), \qquad 1\leq |\a|\leq q,
$$
o\`u $S_{|\a|}(t)$ ne d\'epend que de la longueur de $\a$. Il en r\'esulte 
que les composantes $\xi_\a^{(\mu)}(0)$ d'une d\'eriv\'ee $\xi^{(\mu)}(0)$ d'ordre $\mu$ fix\'e sont de la forme
$\xi_\a^{(\mu)}(0) = c_{|\a|} \s^\a$.

\sk
Pour $\mu\geq 1$ fix\'e, on peut \'ecrire $\xi^{(\mu)}(0) = \sum_{k=1}^q c_k v_k$,
o\`u $v_k$ a pour composantes $v_{k,\a}=\delta_{k |\a|} \, \s^\a$. La courbe param\'etr\'ee $\xi(t)$ engendre 
donc un espace de dimension $q$. Le germe $X$ a la propri\'et\'e $(P)$. 

\bk
R\'eciproquement, si le germe $X$ a la propri\'et\'e (P), le Lemme \ref{conclusion}
permet de supposer qu'il est param\'etr\'e par 
$$
x_\a  = s^\a + \sum_{\mu=q+3}^{+\infty} Q_{\mu,\a}(s), \qquad 1\leq |\a| \leq q,
$$
o\`u les polyn\^omes $Q_{\mu,\a}(s)$ sont homog\`enes de degr\'e $\mu$ et $Q_{\mu,\a}(s)\equiv 0$
si $|\a|=1$, et que les  courbes distingu\'ees sont donn\'ees par  
$$
\xi_\a(t)  = \s^\a t^{|\a|} + \sum_{\mu=q+3}^{+\infty} Q_{\mu,\a}(\s)t^\mu, \qquad 1\leq |\a| \leq q.
$$
Par hypoth\`ese, pour tout $\mu\geq 1$, le vecteur $\xi^{(\mu)}(0)$ appartient \`a l'espace de 
dimension $q$ engendr\'e par $\,\xi'(0),\xi^{(2)}(0), \ldots,\xi^{(q)}(0)$.

\sk
On en d\'eduit l'existence de fonctions $R_{\mu,k}(s)$ telles que 
$$
Q_{\mu,\a}(s) = s^\a R_{\mu,|\a|}(s), \qquad 1\leq |\a|\leq q.
$$
Pour $|\a|=k$ fix\'e, on obtient en particulier $Q_{\mu,(k,0,0,\ldots,0)}(s)s_2^k = Q_{\mu,(0,k,0,\ldots,0)}(s)s_1^k$.
On en d\'eduit que $Q_{\mu,(k,0,\ldots,0)}(s)$ est divisible par $s_1^k$ : $R_{\mu,k}(s)$ est un polyn\^ome.
Le lemme est d\'emontr\'e.
\epf
\bre
\label{comment}
Le r\'esultat pr\'ec\'edent implique en particulier que le germe $X$ est contenu 
dans \og beaucoup \fg\, d'hypersurfaces alg\'ebriques de $\P^N$, 
qu'on obtient en imposant $R_1(s)\equiv 0$, ce qu'on peut faire (alors les param\`etres $s_1,\ldots,s_n$
s'identifient aux coordonn\'ees $x_\a$ avec $|\a|=1$), et en \'eliminant les fonctions $R_k(s)$, $\, k=2,\ldots,q$
dans les \'equations param\'etriques du germe $X$.
Ceci sugg\`ere que, dans le Th\'eor\`eme \ref{Th0} de Bompiani, on pourrait remplacer
les courbes {\em rationnelles normales de degr\'e $q$} par des 
courbes, peu-\^etre alg\'ebriques, {\em lisses et engendrant des espaces de dimension $\leq q$}.
Nous n'avons pas r\'eussi \`a d\'emontrer que c'\'etait le cas.
\ere

\ssct{D\'emonstration du Th\'eor\`eme \ref{Th1bis}}

Soit $X$ un germe de vari\'et\'e lisse en $0\in  \P^N$, v\'erifiant les hypoth\`eses 
du Th\'eor\`eme \ref{Th1bis}.
Il s'agit de montrer qu'il v\'erifie celles du Th\'eor\`eme \ref{Th1}. 
La d\'emonstration est standard.
On se contente 
ici d'en esquisser une en renvoyant \`a \cite{PT1}, \cite{PT2},
o\`u certains arguments sont d\'evelopp\'es dans un cadre plus g\'en\'eral.

\sk
Soit $\text{CR}_q(\P^N)$ la vari\'et\'e de Chow des $1$-cycles rationnels connexes de degr\'e $q$
de $\P^N$, $\;\text{CRN}_q(\P^N)$ l'ouvert de $\text{CR}_q(\P^N)$ des courbes rationnelles
normales de degr\'e $q$ et $\;\text{CRN}_q(\P^N,0)$ l'ensemble des \'el\'ements de $\text{CRN}_q(\P^N)$
qui passent par $0$.

\sk
Soit $\Si$ l'ensemble des courbes 
$C\in \text{CRN}_q(\P^N,0)$ dont le germe 
en $0$ est contenu dans $X$. On montre facilement, voir \cite{PT2},
que $\Si$ est analytique (mais non n\'ecessairement ferm\'ee).
Si $C_\nu\in \Si$  est une suite qui  converge vers 
$C$ dans $\text{CR}_q(\P^N)$, il existe une composante 
irr\'eductible de $C$ qui passe par $0$, dont le germe 
en $0$ est contenu dans $X$. Elle est n\'ecessairement $q$-r\'eguli\`ere, voir \cite{PT1}, 
donc de degr\'e $\geq q$. Il en r\'esulte que 
$C$ est irr\'eductible et appartient \`a $\text{CRN}_q(\P^N)$.
L'ensemble $\Si$ est donc ferm\'e. C'est un sous-ensemble alg\'ebrique de $\text{CR}_q(\P^N)$,
contenu dans $\text{CRN}_q(\P^N,0)$.

\sk
Soit $\Si'$ une composante irr\'eductible de $\Si$,  $\;I' = \{(C,x), \;\; C\in \Si', \;\; x\in C\}$
la vari\'et\'e d'incidence et $\,p: \, I'\rightarrow \P^N$ la projection canonique. 
Par hypoth\`ese on  peut choisir $\Si'$ tel que $X'=p(I')$ contient un germe de dimension $n$. 
La projection $p$ est donc de rang $\geq n$ au point g\'en\'erique de $I'$.

\sk
Si $C_0\in \Si'$, il existe un voisinage $U$ de $C_0$ dans 
$\text{CRN}_q(\P^N,0)$ et une application analytique $\,f: \, U\times \P^1 \rightarrow \P^N$
telle que $\tau\mapsto f(C,\tau)$ est un isomorphisme se $\P^1$
sur $C$ et $f(C,0)=0$. Soit $\tilde{f}$ la restriction de $f$ \`a $(\Si'\cap U)\times \P^1$. 
Si $C\in \Si'\cap U$, par hypoth\`ese, $\tilde{f}$ prend ses valeurs dans $X$
au voisinage de $(C,0)$. Ainsi le rang de $\tilde{f}$ est $\leq n$. Il est donc \'egal \`a $n$ 
au point g\'en\'erique de de ce voisinage. On en d\'eduit que la vari\'et\'e alg\'ebrique 
$X'$ est de dimension $n$ et contient le germe $X$.

\sk
Finalement, soit $v\in \P T_0X$ et $x_\nu\in X\bck\{0\}$ une suite de points qui tend 
vers $0$ dans la direction $v$. Il existe une suite $C_\nu\in \Si'$ telle que 
$x_\nu\in C_\nu$. Quitte \`a extraire une sous-suite de $(x_\nu)$, on peut supposer
que $C_\nu$ converge dans $\text{CRN}_q(\P^N,0)$ vers une courbe $C_0\in \Si'$.
En utilisant une param\'etriqation $f$ comme ci-dessus,
on v\'erifie que la courbe $C_0$ est tangente \`a la direction $v$ en $0$. Le germe $X$
v\'erifie les hypoth\`eses du Th\'eor\`eme \ref{Th1}.

\sct{Une r\'ecurrence}

\ssct{Forme r\'eduite \`a l'ordre $r$}
Les notations g\'en\'erales sont celles du \S 2.1. 
On se donne un germe $X$ de vari\'et\'e lisse de dimension $n$ en $0\in\C^N$.
Dans la d\'efinition suivante, l'espace tangent $T_0X$ est identifi\'e \`a $\C^n$.
\bd
\label{reduit}
Soit $r\in \{1,\ldots,q+2\}$. Le germe $X$  est 
param\'etr\'e sous {\em forme r\'eduite \`a l'ordre $r$} s'il est donn\'e par des \'equations de la forme 
\begin{eqnarray*}
x_\a & = & s^\a + g_\a (s) = s^\a + O(|s|^{q+1}),  \qquad |\a| + r \leq q+1, \\
x_\a & = & s^\a + g_\a (s) = s^\a + O(|s|^{q+2}),  \qquad |\a| + r =    q+2, \\
x_\a & = & s^\a + g_\a (s) = s^\a + O(|s|^{q+3}),  \qquad |\a| + r \geq q+3,
\end{eqnarray*}
bien s\^ur pour $1\leq |\a|\leq q$, et si toute courbe  distingu\'ee tangente \`a $\s\in \C^n\bck\{0\}$ 
est l'image par la param\'etrisation d'un germe de courbe en $0\in \C^n$ 
donn\'e par des \'equations de la forme 
$$
s_i:=\s_i t + u_i(\s,t) = \s_i t+ O(t^{r+1}), \qquad i=1,\ldots,n.
$$
Il est r\'eduit \`a l'ordre $r$ s'il admet une telle param\'etrisation
et r\'eductible \`a l'ordre $r$  s'il existe $\phi\in \text{Aut}\, \P^N$
avec $\phi(0)=0$, tel que $\phi(X)$ soit r\'eduit \`a l'ordre $r$.
\ed
Dans les calculs qui vont suivre, on fera jouer un r\^ole de pivot \`a 
certaines des composantes $x_\a$ de $x\in \C^N$. On introduit 
quelques abbr\'eviations pour les multiindices $\a\in \N^n$ correspondants. 
Pour $\k\in \{1,\ldots,q\}$ et $j\in\{2,\ldots,n\}$, on note :
\begin{eqnarray}
\label{not1}
(\k )       &:= &  (\k,0,\ldots,0), \\  
\label{not2}
(\k-1;j)    &:= &  (\k-1,0,\ldots,0,1,0\ldots,0) \;\;\; \text{o\`u $1$ est \`a la $j$-i\`eme place.}
\end{eqnarray}

\ssct{La r\'ecurrence}

Un germe de vari\'et\'e lisse et $q$-r\'eguli\`ere $X$ en $0\in \C^N$ admet une param\'etrisation de la forme 
$$
x_\a = p_\a(s) + O(|s|^{q+1}), \qquad 1\leq |\a|\leq q,
$$
o\`u les $p_\a(s)$ forment une base de l'espace des polyn\^omes 
nuls en $0$ et de degr\'e $\leq q$. Il existe une  transformation 
lin\'eaire $\phi$ de $\C^N$ telle que le germe $\phi(X)$ soit 
donn\'e par 
$$
x_\a = s^\a + O(|s|^{q+1}), \qquad 1\leq |\a|\leq q.
$$
Autrement dit, $X$ est r\'eductible \`a l'ordre $1$.

\sk
Soit $r\in \{1,\ldots,q+1\}$ et supposons $X$ param\'etr\'e sous forme 
r\'eduite  \`a l'ordre $r$. On \'ecrit sous forme pr\'ecis\'ee 
les composantes $x_\a$ qui ne sont pas {\em a priori} \'ecrites sous forme r\'eduite \`a l'ordre $r+1$,
ainsi que les \'equations des courbes distingu\'ees :
\beq
\label{r1}
x_\a = s^\a + P_\a(s) + O(|s|^{q+2}), \;\; \text{ $P_\a$ homog\`ene de degr\'e $q+1$}, \;\; |\a|+r=q+1,
\eeq
\beq
\label{r2}
x_\a = s^\a + P_\a(s) + O(|s|^{q+3}), \;\; \text{ $P_\a$ homog\`ene de degr\'e $q+2$}, \;\; |\a|+r=q+2,
\eeq
\beq
\label{r3}
s_i=\s_i t + a_i(\s)t^{r+1} + O(t^{r+2}), \qquad i=1,\ldots,n.
\eeq
Le groupe d'\'equations (\ref{r1}) est absent si $r=q+1$, le groupe (\ref{r2}) est absent si $r=1$.

\sk
Sans faire l'hypoth\`ese (P), on peut normaliser la param\'etrisation :
\ble
\label{A}
Si $X$ est r\'eductible \`a l'ordre $r\leq q$, une homographie 
et un changement de param\'etrisation 
permettent de supposer que, par rapport \`a $s_1$ : dans (\ref{r1}) $\;P_{(|\a|)}(s)$ et $P_{(|\a|-1;\,j)}(s)$ 
sont de degr\'es $\leq |\a|-2$ et, dans (\ref{r2}) $\;P_{(|\a|)}(s)$ est de degr\'e $\leq |\a|-1$.
\ele
Alors, si le germe a la propri\'et\'e (P), la param\'etrisation est r\'eduite 
\`a l'ordre $r+1$ :
\ble
\label{B}
On suppose $X$ param\'etr\'e sous forme r\'eduite \`a l'ordre $r\leq q$ et 
normalis\'e comme dans le Lemme \ref{A}. Si $X$ a la 
propri\'et\'e (P), 
tous les  polyn\^omes $P_\a(s)$  de (\ref{r1})--(\ref{r2}) et, apr\`es changement de param\`etre $t$
toutes les fonctions $a_i(\s)$ de (\ref{r3}), sont nuls.
\ele
Enfin :
\ble
\label{C}
Si $X$ est r\'eduit \`a  l'ordre $q+1$ et a la propri\'et\'e (P), il admet une param\'etrisation 
de la forme (\ref{normal}). Les courbes distingu\'ees de $X$ sont alors les images 
par cette param\'etrisation des germes de droites en $0\in \C^n$.
\ele
Le Lemme \ref{conclusion} en d\'ecoule. Si le germe $X$ a la propri\'et\'e (P),
les deux premiers lemmes montrent, par r\'ecurrence sur $r$, qu'il est
r\'eductible \`a l'ordre $q+1$. Le dernier lemme donne alors le r\'esultat.

\sk
Reste \`a d\'emontrer les trois lemmes pr\'ec\'edents.

\ssct{D\'emonstration du Lemme \ref{A}}

Soit $r\in \{1,\ldots,q\}$. On suppose $X$ param\'etr\'e sous forme 
r\'eduite  \`a l'ordre $r$, la propri\'et\'e $(P)$ n'intervient pas.
La normalisation est obtenue gr\^ace aux transformations suivantes.

\sk
1)  Une homographie $\phi$ correspondant \`a un changement de coordonn\'ees de la forme  
$$
x_\a = x'_\a(1+\wh{G}(x'))^{-1}, \;\; \;1\leq |\a|\leq q, \;\;\;\;  \text{o\`u} \;\;\;\;  \wh{G}(x')= \sum_{|\b|=r} c_\b x'_\b
$$
est {\em une forme lin\'eaire de poids $r$}. On lui associe le polyn\^ome 
homog\`ene de degr\'e $r$ $\;G(s) = \sum_{|\b|=r} c_\b s^\b$. L'application $\wh{G}\mapsto G$ est une bijection 
de l'espace des formes lin\'eaires de poids $r$ sur l'espace des 
polyn\^omes homog\`enes de degr\'e $r$.

\sk
2) Un changement de param\`etre de la forme $\;s_j = s'_j + H_j(s')$, $\;j=1,\ldots,n$,
o\`u les $H_j(s')$ sont des polyn\^omes homog\`enes de degr\'e $r+1$ ;

\sk
3) Une transformation lin\'eaire anodine dont l'objet est 
de mettre \`a nouveau sous forme r\'eduite les composantes $x_\a$ dont le  poids $|\a|$ v\'erifie $|\a|+r\leq q$.

\bk
La param\'etrisation initiale de $X$ est de la forme 
$x_\a = s^\a + O(|s|^{m_{|\a|}})$, o\`u l'entier 
$m_{|\a|}\in \{q+1,q+2,q+3\}$
est au moins \'egal \`a $|\a|+1$.
Sur l'image $\phi(X)$ de $X$ par une homographie 1), on a
$$
x_\a(1+\wh{G}(x))^{-1} = s^\a + O(|s|^{m_{|\a|}}), \qquad 1\leq |\a| \leq q,
$$
donc en particulier $x_\a = s^\a + O(|s|^{|\a|+1})$ puis $\wh{G}(x)=G(s)+O(|s|^{r+1})$.
On en d\'eduit que la param\'etrisation de $\phi(X)$ induite par la param\'etrisation
initiale est de la forme 
$$
x_\a = s^\a + O(|s|^{|\a|+r}) + O(|s|^{m_{|\a|}}), \qquad 1\leq |\a|\leq q.
$$
Un changement de variable 2) ne change 
rien \`a cette forme ni \`a celle des \'equations des courbes distingu\'ees.
Ci-dessus, une composante
$x_\a$  est \'ecrite sous forme r\'eduite \`a l'ordre $r$ si 
$|\a|+r\geq m_{|\a|}$, autrement dit si $|\a|+r\geq q+1$.
Toutefois, une transformation lin\'eaire  
$$
x_\a = l_\a(x'), \;\; |\a|+r\leq q, \qquad x_\a=x'_\a, \;\; |\a|+r\geq q+1,
$$
r\'eduit \`a nouveau \`a l'ordre $r$ les  composantes $x_\a$ dont le poids $|\a|$
v\'erifie $|\a|+r\leq q$, sans modifier les composantes dont le poids $|\a|$ 
v\'erifie $|\a|+r\geq q+1$.

\sk
Ces consid\'erations montrent qu'il suffit  d'analyser l'influence d'une homographie 1) et d'un changement de
param\`etre 2) sur les composantes (\ref{r1})--(\ref{r2}), c'est-\`a-dire 
celles dont le poids $|\a|$ v\'erifie $|\a|+r\in \{q+1,q+2)$.
Par hypoth\`ese, 
$$
x_\a = s^\a + P_\a(s) + O(|s|^{|\a|+r+1}), \qquad |\a|+r\in \{q+1,q+2\},
$$
o\`u les polyn\^omes $P_\a(s)$ sont homog\`enes de degr\'e $|\a|+r$.

\sk
Si $\phi$ est l'homographie 1), les composantes analogues de $x\in \phi(X)$ v\'erifient 
les relations $x_\a (1+\wh{G}(x))^{-1} = s^\a + P_\a(s) + O(|s|^{|\a|+r+1})$, 
donc sont donn\'ees par 
$$
x_\a = s^\a + P_\a(s) + s^\a G(s) +  O(|s|^{|\a|+r+1}).
$$
Un changement de param\`etre 2) les met ensuite sous la forme 
\beqn
x_\a      & = &  s^\a + P'_\a(s) + O(|s|^{|\a|+r+1}), \\
P'_\a(s)  & = &  P_\a(s) + s^\a G(s) + \sum_{j=1}^n  \fr{\a_js^\a}{s_j}H_j(s), \qquad |\a|+r\in \{q+1,q+2\}.
\eeqn

\sk
Pour $|\a|+r=q+1$, on obtient en particulier le syst\`eme d'\'equations homog\`enes :
\beqn
P'_{(|\a|)}(s) 
& = & 
P_{(|\a|)}(s) + s_1^{|\a|} G(s) + |\a| s_1^{|\a|-1}H_1(s), \\
P'_{(|\a|-1;\,j)}(s) 
& = & 
P_{(|\a|-1;\,j)}(s) + s_1^{|\a|-1}s_jG(s) + (|\a|-1)s_1^{|\a|-2}s_jH_1(s) + s_1^{|\a|-1}H_j(s).
\eeqn
Quel  que soit $G(s)$, homog\`ene de degr\'e $r$,
on peux choisir 
$H_1(s),\ldots,H_n(s)$, homog\`enes de degr\'e $r+1$, tels que $P'_{(|\a|)}(s)$ et les $P'_{(|\a|-1;\,j)}(s)$
soient de degr\'es $\leq |\a|-2$ par rapport \`a $s_1$. En particulier, 
$H_1(s) = - s_1G(s))/|\a| + A(s)$ o\`u $A(s)$ d\'epend des donn\'ees. 
Si $r=1$, on peut choisir $G(s)\equiv 0$ et la normalisation est termin\'ee.
On suppose $r\geq 2$.

Pour $|\a|+r=q+2$, on obtient 
$P'_{(|\a|)}(s)  = P_{(|\a|)}(s) + s_1^{|\a|} G(s) + |\a|s_1^{|\a|-1} H_1(s)$, donc :
$$
P'_{(|\a|)}(s)  = P_{(|\a|)}(s) - s_1^{|\a|} G(s)/(|\a|-1) + B(s),
$$
o\`u $B(s)$ d\'epend des donn\'ees. On peut choisir (il y a une seule  solution) $G(s)$, 
homog\`ene de degr\'e $r$, de telle sorte que $P'_{(|\a|)}(s)$
soit de degr\'e $\leq |\a|-1$  par rapport \`a $s_1$. Ceci ach\`eve la normalisation 
annonc\'ee dans le Lemme \ref{A}.

\ssct{Traduction de la propri\'et\'e $(P)$}

Soit $X$ un germe de vari\'et\'e lisse $q$-r\'eguli\`ere en $0\in \C^N$, 
param\'etr\'e sous  forme r\'eduite \`a l'ordre $1$ : $\;x_\a = s^\a + O(|s|^{q+1})$
pour $\;1\leq |\a|\leq q$.
Il est clair que, si $\xi(t) = (\xi_\a(t))_{1\leq |\a|\leq q}$ avec $\xi(0)=0$ et $\xi'(0)\neq 0$,
est un germe de courbe param\'etr\'ee contenue dans $X$, alors  
$\xi'(0),\ldots,\xi^{(q)}(0)$ sont lin\'eairement ind\'ependants. 
Si ce germe est une courbe distingu\'ee, il est donc contenu 
dans l'espace engendr\'e par ces d\'eriv\'ees et pour tout $\mu\geq 1$ donn\'e, 
il existe $c_1,\ldots,c_q\in \C$ tels que :
$$
\xi^{(\mu)}(0) = \sum_{k=1}^q c_k \, \xi^{(k)}(0).
$$
On se donne un entier $\k\in \{1,\ldots,q\}$ et un entier $r\geq 1$.
On suppose que les composantes de poids $|\a|=\k$ sont donn\'ees par 
\beq
\label{h1}
x_\a = s^\a + P_\a(s) + O(|s|^{\k+r+1}), \qquad |\a|=\k,
\eeq
o\`u les $P_\a(s)$ sont homog\`enes de degr\'e $\k+r$, 
et que les courbes distingu\'ees sont les images par la param\'etrisation
de germes donn\'es par des \'equations de la forme :
\beq
\label{h2}
s_i=\s_i t + a_i(\s)t^{r+1} + O(t^{r+2}), \;\; i=1,\ldots,n.
\eeq
De (\ref{h1}) et (\ref{h2}), on tire les composantes $\xi_\a(t)$ de poids $|\a|=\k$
d'une courbe distingu\'ee :
$$
\xi_\a(t)  = 
\s^\a  t^\k + \left( \sum_{i=1}^n a_i(\s) \, \fr{\a_i\s^\a}{\s_i} + P_\a(\s)\right) t^{\k+r} + O(t^{\k+r+1}).
$$
Si $\k+r\geq q+1$, la seule d\'eriv\'ee non nulle $\xi_\a^{(k)}(0)$ d'ordre $k\leq q$ est $\xi_\a^{(\k)}(0) = \k! \, \s^\a$. 
De $\xi^{(\k+r)}(0) = \sum_{k=1}^q c_k \, \xi^{(k)}(0)$,
on tire la traduction suivante de la condition (P).

\sk
{\em Si le germe $X$ est param\'etr\'e sous forme r\'eduite \`a l'ordre $1$ et a la propri\'et\'e (P),
si les composantes $x_\a$ de poids $|\a|=\k$ de $x\in X$ sont donn\'ees par (\ref{h1}) 
et les courbes distingu\'ees par (\ref{h2}), il existe un scalaire $c(\s)=c_{\k,r}(\s)$ tel que :}
\beq
\label{f1}
\k+r \geq q+1 \; \Rightarrow \; \sum_{i=1}^n a_i(\s)  \fr{\a_i \s^\a}{\s_i} + P_\a(\s) = c(\s) \s^\a, \qquad |\a|=\k.
\eeq

\ssct{D\'emonstration du Lemme \ref{C}}

On suppose que  $X$ est param\'etr\'e sous forme r\'eduite \`a l'ordre $q+1$ : 
$$
x_\a = s^\a + O(|s|^{q+2}), \;\;  |\a|=1, \qquad x_\a = s^\a + O(|s|^{q+3}), \;\;  2\leq  |\a| \leq q.
$$
En prenant $(x_\a)_{|\a|=1}$ comme nouveau param\`etre $s=(s_1,\ldots,s_n)$,
on obtient de plus :
$$
x_\a = s^\a, \qquad |\a| = 1.
$$
Les courbes distingu\'ees sont encore donn\'ees par (\ref{h2}) avec $r = q+1$.
Si par exemple $\s_1\neq 0$, on prend $s_1/\s_1$ comme 
nouveau param\`etre $t$. Alors $a_1(\s)\equiv 0$.

\sk
On suppose les courbes distingu\'ees donn\'ees par (\ref{h2}) avec $r \geq q+1$
donn\'e quelconque. Si $|\a|=1$, les relations (\ref{f1}) donnent $a_i(\s)=c(\s)\s_i$
puisque $P_\a(s)\equiv 0$, donc $a_i(\s)\equiv 0$, $\, 1\leq i\leq n$.
Par r\'ecurrence sur $r$, on obtient $s(t)\equiv \s t$ et le Lemme \ref{C}.

\ssct{D\'emonstration du Lemme \ref{B}}

On suppose que le germe $X$ est param\'etr\'e sous forme 
r\'eduite  \`a l'ordre $r\leq q$, qu'il est normalis\'e comme dans le Lemme \ref{A}  et qu'il 
a la propri\'et\'e $(P)$. Rappelons les notations. On a :
\beq
\label{p1}
x_\a = s^\a + P_\a(s) + O(|s|^{|\a|+r+1}), \qquad |\a|+r\in \{q+1,q+2\},
\eeq
o\`u les $P_\a(s)$ sont homog\`enes de degr\'e $|\a|+r$. Les courbes distingu\'ees sont donn\'ees par
\beq
\label{p2}
s_i=\s_i t + a_i(\s)t^{r+1} + O(t^{r+2}), \qquad i=1,\ldots,n.
\eeq
Si par exemple $\s_1\neq 0$, on prend $s_1/\s_1$ comme nouveau param\`etre $t$. Alors :  
$$
a_1(\s)\equiv 0.
$$
En notant $s$ au lieu de $\s$, pour $s$ voisin de $(1,0,\ldots,0)$, l'implication (\ref{f1}) donne  : 
\beq
\label{f2}
\sum_{i=2}^n a_i(s) \, \fr{\a_i s^\a}{s_i} + P_\a(s) = c_{|\a|}(s) s^\a, \qquad |\a|+r\in \{q+1,q+2\}.
\eeq
On note $|\a|=\k$ et on traite d'abord simultan\'ement les deux groupes d'\'equations correspondant 
respectivement \`a $\k+r=q+1$ et $\k+r=q+2$. 

\sk
Avec les notations (\ref{not1})--(\ref{not2}), pour $\a=(\k)$, on obtient $c_\k(s)s_1^\k  = P_{(\k)}(s)$
car $a_1(\s)\equiv 0$.
En choisissant $\a=(\k-1;j)$, on obtient $\,a_j(s)s_1^{\k-1} + P_{(\k-1;j)}(s) = c_\k (s) s_1^{\k-1}s_j$
donc 
\beq
\label{cond1}
a_j(s)s_1^\k = P_{(\k)}(s)s_j -  P_{(\k-1; \,j)}(s)s_1, \qquad j=2,\ldots,n.
\eeq
En multipliant  les deux membres de (\ref{f2}) par $s_1^\k$, on obtient  :
\beq
\label{cond2}
P_\a(s)s_1^\k  = (1-\sum_{i=2}^n \a_i)P_{(\k)}(s)s^\a  
+\sum_{i=1}^n P_{(\k-1;\,i)}(s) \fr{\a_i s_1 s^\a}{s_i}, \;\;\; |\a|=\k. 
\eeq
{\em Il s'agit d'une identit\'e entre polyn\^omes :  $ \a_is^\a/s_i$ est un mon\^ome, nul si $\a_i=0$.}

Le mon\^ome $s_1^\k$ divise le second membre de (\ref{cond2}).
Pour $j\geq 2$ donn\'e, en choisissant d'abord $\a$, de longueur $|\a|=\k$, tel que $\a_j=\k$, puis 
tel que $\a_j=\k-1$ et $\a_1=1$, on obtient : 
\begin{eqnarray}
\label{div1}
\text{$s_1^\k\;$  divise} \;\;  & & 
(1-\k)P_{(\k)}(s)s_j  + \k P_{(\k-1;j)}(s)s_1, \;\;\;\;\\
\label{div2}
\;\;\;\;\;\;\text{$s_1^\k\;$ divise} \;\; & & 
(2-\k)P_{(\k)}(s)s_j s_1 + (\k-1)P_{(\k-1;\,j)}(s)s_1^2. 
\end{eqnarray}

\sk
On suppose d'abord $|\a|+r:=\k+r=q+1$. 
Par les hypoth\`eses du Lemme \ref{B}, les polyn\^omes $P_{(\k)}(s)$ et $P_{(\k-1;\,j)}(s)$
sont de degr\'e $\leq \k-2$ par rapport \`a $s_1$.
D'autre part, il  r\'esulte de (\ref{div1})--(\ref{div2}) que 
$s_1^{\k-1}$ divise les deux polyn\^omes 
$$
(1-\k)P_{(\k)}(s)s_j  + \k P_{(\k-1;\,j)}(s)s_1,
\;\;\;
(2-\k) P_{(\k)}(s)s_j + (\k-1)P_{(\k-1;\,j)}(s)s_1,
$$
donc aussi le polyn\^ome $P_{(\k)}(s)$ : il est  nul.
Il r\'esulte alors de (\ref{div1}) que $s_1^{\k-1}$ divise 
les polyn\^omes $P_{(\k-1;\,j)}(s)$ : ils sont nuls. 
Alors (\ref{cond1})--(\ref{cond2}) montre que toutes les fonctions
$a_j(s)$ et tous les polyn\^omes $P_\a(s)$ avec $|\a|+r=q+1$
sont nuls.

\sk
Si $r\neq 1$, on suppose enfin $|\a|+r=q+2$.
Par hypoth\`ese le polyn\^ome $P_{(\k)}(s)$ 
est de degr\'e $\leq \k-1$ par rapport \`a $s_1$.
Comme les fonctions $a_i(\s)$ sont nulles, (\ref{cond1})
donne $P_{(\k-1; \,j)}(s)s_1 = P_{(\k)}(s)s_j$ et on obtient 
que $s_1^\k$ divise le polyn\^ome $P_{(\k)}(s)$. Il est nul 
donc aussi les $P_{(\k-1;\,j)}(s)$ et par (\ref{cond2}) tous les 
$P_\a(s)$. Le Lemme \ref{B} est d\'emontr\'e.

\sct{Appendice : commentaires sur l'article d'E. Bompiani}

\ssct{Le th\'eor\`eme g\'en\'eral}

Si $Y\subset \P^N$ est un germe lisse de vari\'et\'e analytique, on note 
$Y_x(q)$ son espace osculateur \`a l'ordre $q$ en $x\in Y$. 

\sk
Soit $X$ un germe lisse $q$-r\'egulier. Une courbe analytique lisse 
$C\subset X$ est {\em quasi-asymptotique de type $(q-1,q+1)$} si
pour tout $x\in C$, l'espace $C_x(q+1)$ est contenu dans l'espace 
engendr\'e par $X_x(q-1)$ et $C_x(q)$.
Le r\'esultat principal de \cite{Bo} est le suivant :
\bt[Bompiani]
Soit $X$ un germe de vari\'et\'e lisse $q$-r\'egulier de dimension $n$.
Si $X$ contient $\infty^{2n-2}$ courbes quasi-asymptotiques de type
$(q-1,q+1)$,  $\;X$ est une vari\'et\'e de Veronese d'ordre $q$ et ces 
courbes sont rationnelles normales de degr\'e $q$.
\et
Nous ne red\'emontrons pas ce r\'esultat. 
La presque totalit\'e de \cite{Bo} est consacr\'ee \`a un calcul formel
difficile dont l'auteur d\'eduit que, sous les hypoth\`eses du th\'eor\`eme : 

\sk
{\em toute courbe quasi-asymptotique de type
$(q-1,q+1)$ contenue dans $X$ engendre un espace de dimension $q$ ;
de plus deux telles courbes \og infiniment voisines \fg\, engendrent 
un espace de dimension $2q$.}

\sk
Le calcul est d\'etaill\'e d'abord dans le cas o\`u $X$ est une surface,
puis dans le cas $n=3$ pour $q=2$ et de fa\c{c}on plus succincte 
pour $q$ quelconque.
Le cas g\'en\'eral est laiss\'e au lecteur \footnote{
{\em ``\`E certo che il procedimento dimostrativo rimane inaltero,
e quindi la verifica dell'asserto presenter\`a, al pi\`u, difficolt\`a
materiali inessenziali: per questo mi \`e parso inutile 
farla''}, \cite{Bo}~page 474. }.
Ce calcul est par exemple forc\'ement plus d\'elicat que 
celui qu'on a fait pour d\'emontrer le Th\'eor\`eme \ref{Th1}.
En effet, pour d\'emontrer le th\'eor\`eme de Bompiani,
il faut utiliser non seulement l'existence de 
courbes quasi-asymptotiques au point de base du calcul, 
mais aussi l'existence de courbes quasi-asymptotiques 
qui passent pr\`es du point de base. La difficult\'e 
est alors de d\'eterminer les \'equations de compatibilit\'e 
utiles que cela impose aux \'equations d'abord obtenues.

\bk
La deuxi\`eme \'etape de la d\'emonstration est courte et g\'eom\'etrique 
(voir \cite{Bo}, pages 462--463 pour $n=2$ et pages 473--474 pour $n=3$).
\`A partir du r\'esultat d\'ej\`a obtenu, l'auteur montre que,
sous les hypoth\`eses du th\'eor\`eme :

\sk
{\em  les courbes quasi-asymptotiques de type
$(q-1,q+1)$ contenues dans $X$ sont des courbes rationnelles normales 
de degr\'e $q$.}

\bk
L'auteur est ainsi ramen\'e \`a d\'emontrer le Th\'eor\`eme \ref{Th0}.
Si $X$ est une surface, voir \cite{Bo} page 464 : et puisque la surface 
contient un faisceau
(donc une infinit\'e) rationnel de courbes rationnelles de degr\'e $q$,
elle est certainement rationnelle et comme elle engendre $\P^{q(q+3)/2}$,
il est bien connu que c'est la surface repr\'esentative de toutes les courbes 
planes de degr\'e $q$. L'\'enonc\'e est donc connu dans ce cas, peut-\^etre aussi, 
la phrase pr\'ec\'edente est ambigu\"{e}, le Th\'eor\`emes \ref{Th1} et le Th\'eor\`eme \ref{Th1bis}.

\ssct{La d\'emonstration du Th\'eor\`eme \ref{Th0} dans \cite{Bo}}

Hormis le cas connu $n=2$, l'auteur donne page 470 une d\'emonstration \og synth\'etique \fg\, 
du Th\'eor\`eme \ref{Th0} pour $n=3$, $q=2$. Le cas $q\geq 3$ fait l'objet 
d'une ligne de la page 475. Trouvait-il le r\'esultat trop banal ou trop \'evident ?

\sk
Je pense que sa d\'emonstration dans le cas $q=2$ pr\'esente une lacune
\footnote{Si par chance ces lignes ont un lecteur et si 
l'article de Bompiani lui est familier, je lui serais tr\`es reconnaissant 
de me donner son opinion sur ce sujet.}.
Je reprends ci-dessous une partie de l'argumentation  en mettant en \'evidence 
le point qui m'intrigue : une certaine propri\'et\'e,
utilis\'ee sans commentaire, exigerait selon moi
une d\'emonstration. 

\sk
La strat\'egie 
pour d\'emontrer le r\'esultat est de construire un syst\`eme de 
diviseurs de $X$ dont l'\'el\'ement g\'en\'erique 
est une vari\'et\'e de Veronese de dimension $n-1$ et d'ordre $q$,
tel que par tout $n$-uplet g\'en\'erique de points de $X$
passe un diviseur du syst\`eme.
Ce syst\`eme est un syst\`eme lin\'eaire qui r\'esout le probl\`eme :
l'application $X\dasharrow \P^N$ associ\'ee envoie les courbes 
rationnelles normales de degr\'e $q$ de $X$ sur les droites de $\P^N$.
Par exemple, dans le cas $n=2$, il suffit de prendre le syst\`eme de 
ces courbes. Dans la suite, je consid\`ere seulement la partie de 
la d\'emonstration qui concerne la construction d'un tel syst\`eme.

\ssct{Le point cl\'e si $n=3$, $q=2$}

Voici comment je comprends l'argument essentiel de la d\'emonstration
de \cite{Bo}. On part d'une vari\'et\'e $X$ de dimension $3$ 
qui engendre un $\P^9$, telle que par toute paire $(x,y)\in X^2$
g\'en\'erique passe une conique propre $C(x,y)\subset X$. Il s'agit de montrer que 
par trois points de $X$ passe une surface de Veronese d'ordre $2$
(qui engendre donc un $\P^5$) contenue dans $X$.

\sk
Soit $(a,b,c)$ un triplet g\'en\'erique de points de $X$ et $S=S(a,b,c)$
la surface d\'ecrite par les coniques de $X$ qui passent par $a$ 
et rencontrent la conique $C(b,c)$. Le plan tangent $S_a(1)$
est engendr\'e par les droites tangentes $C(a,b)_a(1)$ et $C(a,c)_a(1)$.
Quand $x$ d\'ecrit $C(b,c)$, le plan de la conique $C(a,x)$
est engendr\'e par sa tangente en $a$ et le point $x$. Il en r\'esulte 
que la surface $S$ est contenue dans le $\P^5$ engendr\'e par le plan $S_a(1)$
et le plan de $C(b,c)$.

\sk
Si $x\in C(b,c)$, on peut appliquer le m\^eme raisonnement  
\`a $S'=S(x,a,c)$. On a $S'_x(1)=S_x(1)$
puisque les deux espaces contiennent $C(x,a)_x(1)$ et $C(b,c)_x(1)$. 
Pour tout $y\in C(a,c)$, la conique $C(x,y)$ est contenue 
{\em dans le m\^eme $\P^5$ qui contient $S$}, puisque ce $\P^5$
est aussi engendr\'e
par $S_x(1)=S'_x(1)$ et le plan de $C(a,c)$.
Finalement, les coniques de $X$ qui rencontrent 
\`a la fois $C(a,c)$ et $C(b,c)$ engendrent un $\P^5$
et d\'ecrivent une surface qui contient $\infty^2$ coniques, une surface de Veronese
d'ordre $2$.

\sk
Comme Bompiani, nous avons introduit sans commentaire l'espace tangent 
en $a$ \`a la surface $S(a,b,c)$,
mais je ne vois pas pourquoi, compte tenu de la fa\c{c}on dont elle est d\'efinie,
la r\'eponse \`a la question suivante serait positive.

\sk
{\bf \em La surface $\bf S(a,b,c)$ est-elle lisse au point $\bf a$ ?}

\ssct{Le cas g\'en\'eral}

Bompiani ne donne pas d'indication sur sa fa\c{c}on de traiter le cas
g\'en\'eral. Connaissait-il le r\'esultat par ailleurs ?
Voici une esquisse de \og d\'emonstration \fg\, qui repose sur la m\^eme 
id\'ee et le m\^eme genre 
de p\'etition de principe que la pr\'ec\'edente. C'est  peut-\^etre la d\'emonstration
qu'il avait en t\^ete.

\sk
Soit $X\subset \P^N$ une vari\'et\'e de dimension $n$, $\,q$-r\'eguli\`ere 
au point g\'en\'erique, telle que pour tout $(x,y)\in X^2$ g\'en\'erique, il existe une
courbe rationnelle normale $C(x,y)\subset X$ de degr\'e $q$.
\`A tout $n$-uplet g\'en\'erique $(x_1,\ldots,x_n)\in X^n$ on associe les sous-vari\'et\'es
$$
Y(x_1)\subset Y(x_1,x_2) \subset \cdots \subset Y(x_1,\ldots,x_{n-1}) \subset Y(x_1,\ldots,x_n),
$$
o\`u $Y(x_1)$ est le point $x_1$, $Y(x_1,x_2)$ la courbe $C(x_1,x_2)$ 
et plus plus g\'en\'eralement, $Y(x_1,\ldots,x_{r+1})$ la vari\'et\'e de dimension
$r$ d\'ecrite ou presque par les courbes $C(x_{r+1},z)$ quand $z$ d\'ecrit $Y(x_1,\ldots,x_r)$.
Il s'agit de montrer que, pour $(x_1,\ldots,x_n)\in X^n$ g\'en\'erique,
l'hypersurface $Y=Y(x_1,\ldots,x_n)$ 
est une vari\'et\'e de Veronese d'ordre $q$. 

\sk
{\bf \em On fait l'hypoth\`ese : pour $\bf (x_1,\ldots,x_n)$ g\'en\'erique,  $\bf Y(x_1,\ldots,x_n)$
 est lisse au point $\bf x_n$ (\og point de base \fg\, de la construction).}

\sk
Tout est alors facile. On peut supposer que $X$ est lisse et $q$-r\'egulier 
en $x_n$ et que $Y$ est lisse donc $q$-r\'egulier en $x_n$.
Comme la vari\'et\'e $Y$ est recouverte ou presque par les courbes $C(x_n,y)$, $\,y\in Y$,
elle est  contenue dans son osculateur $Y_{x_n}(q)$. 
Comme elle est $q$-r\'eguli\`ere en $x_n$, elle engendre un espace 
$\lan \, Y\,  \ran$ de dimension ${n-1+q \choose q} - 1$.

\sk
Supposons qu'une courbe g\'en\'erique $C(x,x')$ rencontre $Y$ en au moins deux points.
Soit $x\in X\bck Y$ g\'en\'erique tel que $X$ est lisse et $q$-r\'egulier en $x$.
La vari\'et\'e $X$ est recouverte ou presque par les courbes $C(x,y)$ qui passent par $x$
et un point $y\in Y$. Celles-ci ont deux points dans $Y$ donc 
l'espace $\P^N$ qu'elles engendrent est contenu
dans l'espace engendr\'e par $\lan \, Y\, \ran$ et $X_x(q-2)$. On obtient la contradiction :
$$
{ n+q \choose q } \leq  { n-1+q \choose q } + {n+q-2 \choose q-2}.
$$
Une courbe g\'en\'erique $C(x,x')$ rencontre donc une vari\'et\'e $Y=Y(x_1,\ldots,x_n)$ 
g\'en\'erique en un point. Autrement dit, pour $(y_1,y_2)\in Y^2$ g\'en\'erique,
la courbe $C(y_1,y_2)$ qui a deux points dans $Y$ est contenue dans $Y$ : $Y$ est 
une vari\'et\'e de Veronese d'ordre $q$. Des variations sont possibles, 
par exemple 
consid\'erer les vari\'et\'es $Y$ d\'ecrites par les courbes tangentes en un 
point $x\in X$ \`a un hyperplan donn\'e de $T_xX$. On obtiendrait la m\^eme 
conclusion, \`a condition de poser la m\^eme p\'etition de principe :
que l'hypersurface $Y$ est lisse au point $x$.

\end{document}